\documentclass[11pt,english]{article}
\makeatletter
\usepackage{geometry}
\usepackage{float}
\usepackage{graphicx}
\usepackage{amsmath,amsthm}
\usepackage[american]{babel}
\newcommand{\doublespacing}{\let\CS=\@currsize\renewcommand{\baselinestretch}{1.2}\tiny\CS}

\date{}
\usepackage{hyperref}
\hypersetup{colorlinks=true,linkcolor=blue,citecolor=blue}
\begin{document}
\newpage
\doublespacing
\doublespacing
\title {\Large {  Optimum Production for a heaped stock  dependent breakable item through variational principle}}
\author{{J. N. Roul $^{a},$ K. Maity  $^{b}, $   S. Kar $^{c},$ and M. Maiti $^{d}$}\\\\
  $^{a}$Patha Bhavana, Visva-Bharati, Santiniketan-731235, West Bengal, India,
\\   Email$-$jotin2008@rediffmail.com \\\\
  $^{b}$ Department of Mathematics, Mugberia
 Gangadhar Mahavidyalaya, Bhupatinagar,\\   Purba
Medinipur-721425, W.B., India, Email$-$kalipada\_maity@yahoo.co.in\\\\
   $^{c}$ Department of Mathematics, National Institute of
Technology, Durgapur-713209\\ W.B., India, Email$-$kar\_s\_k@yahoo.com\\\\
  $^{d}$Department of Applied Mathematics with Oceanology and Computer
 Programming,\\
 Vidyasagar University, Midnapore-711201,W.B., India,
Email$-$mmaiti2005@yahoo.co.in}

\maketitle
\begin{abstract}
Breakability rate of fragile item depends on the accumulated stress of heaped stock level. So   breakablility rate can be considered as dependent parameter  of stock variable.  The unit production cost is a function of production rate and also dependent on raw material cost, development cost  and wear-tear cost. The holding cost is assumed to be non-linear, dependent on time. Here optimal control problem  for a fragile  item under finite  time horizon is considered.  The profit  function which consists of  revenue, production  and holding costs is formulated as a Fixed-Final Time and Fixed State System(cf. Naidu (2000)[\ref{nai2000}]) optimal control problem with finite time horizon. Here  production rate is unknown and considered as a control variable and  stock level is taken as a state variable. It  is formulated to optimize the production rate so that total profit is  maximum.  As particular cases,  models are evaluated with and without  breakability. The models are  solved by using conventional Variational Principle  along with the non-linear optimization technique-Generalised Reduced Gradient Method (LINGO 12.0).  The optimum results are illustrated  both  numerically and graphically. Some sensitivity analysis on breakability coefficient are presented.
\end{abstract}
\section{Introduction}

Available items  in the market can be broadly classified into two categories-damageable items and non-damageable items. Damageable items again can be classified into two sub-categories -deteriorating and breakable items. Deteriorating items are deteriorated with time.  Normally, seasonal goods like fruits, vegetables, X-mas cake, etc., are deteriorating in nature. Demands of these items normally exist in the market for a finite time and obviously these types of demands are dependent on time.\\

   Fashionable/decorating items are made of glass, ceramic, etc., and demand of these types of items continues over a fixed duration only. Sale of these fashionable products increases with the showing of stock. But manufacturers  face a conflicting situation in their business as breakability increases with the increase of piled stock.   Inventory practitioners  are not interested with this types of items for their breakability. Few articles in this direction have been focused by Mandal and Maiti 2000 [\ref{mandal2000}], Guchhait et al. (2010,2013) [\ref{guchhait2010},\ref{guchhait2013}].\\
Almost all inventory models are formulated with constant  holding cost ( Sarkar et al., 2011 [\ref{sar2011}; Roul, et al. 2017 [\ref{rou2017}] ). In reality, due to inflation, bank interest, rental charges, preservation cost, etc., it increases with time. Thus some factors are contributed to the holding cost change with time (Giri et al., 1996 [\ref{gir1996}], Alfares 2007 [\ref{alfares2007}],  Urban 2008 [\ref{urban2008}] ) and others remain constant. Also set-up cost depends on production rate as high production rate require sophisticated modern machineries. In this paper, holding and set-up costs are considered as functions of time and production rate respectively as a particular case. This is very simple application of breakability  in production inventory problem via variational principle. \\

In this paper, a production-inventory model with imperfect production process is considered for a breakable  item over a finite time horizon. The production rate varies with time.  Set-up cost is partially production rate dependent and holding cost is also partially time dependent. The unit production cost is a function of production rate, raw material cost, labour charge, wear and tear cost. The
models are formulated as optimal control problems for the maximization of total profits over the planning horizon and solved using variational principle. The profit is maximized using the optimization technique-Generalised Reduced Gradient method(LINGO 12.0) (cf. Gabriel and Ragsdell 1977)[\ref{gab1977}]. The models are numerically illustrated and optimum results are presented in tabular form and graphically . Some sensitivity analysis on breakability coefficient are given.
\newpage
\section{Proposed Optimal Control Models }
For  the optimal control problem under finite time horizon, following assumptions and notations are used.

\subsection{Assumptions}
\begin{itemize}
\item It is a single period production inventory model with finite time horizon $t\epsilon[0,T]$;
\item  Break-ability rate is function of stock;
\item Shortages are not allowed;
\item There is no repair or replacement of breakable units over whole time period;
\item Unit production cost depends on produced-quantity, raw material, wear-tear and development costs;
\item  Deterministic  variational principle is  used for the model.
\end{itemize}
\subsection{Notations}
{$T :$} total time for the system;\\
{$u(t):$} production rate at time $t$ which is taken as a control variable;\\
{$x(t):$} stock level at  time t which is a state variable;\\
{$h(t) :$}=$a+b t^n( a\,\,\,\, \mbox{and} \,\,\,\,b \,\,\,\,\mbox{are constants })$ holding cost per unit;\\ 
{$L :$} fixed cost like labour, energy, etc;.\\  
{$N :$} cost of technology, design, complexity, resources, etc;.\\  
{$c_{10}:$} constant material cost;\\
{$d_{1}:$} constant demand at initial stage;\\
{$d(t)=d_{1}+d_2 t+d_3 t^2:$} demand function with $d_2$ and $d_3$ as constants.\\
{$c_d(l)=N+L: $} development cost to improve the quality of the product;\\
{$\beta_{10}:$} wear-tear cost for the system.\\
{$p:$} is the selling price of unit quantity. \\
$ c_{u}(t)=c_{10} +\frac{c_{d}(l)}{u(t)}+\beta_{10} u(t)$ is the unit production cost.

\subsection{Model-1: Formulation of optimal control  models for quadratic  demand with holding cost $h(t)=a+bt^n$}
\indent Let us consider  a production system with production rate $ u(t)$, demand $d(t)$ and stock dependent breakability rate $B(x)$. In this model,  the differential equation for stock level
{$x(t)$} regarding above system
during a fixed time-horizon, T is
\begin{eqnarray} \frac{d}{dt}\bigg(x(t)\bigg)= u (t)-d(t)-B(x)
\end{eqnarray}
where $ d(t)=d_1 +d_2 t+d_3 t^2$ and $B(x)\equiv$ $B(x(t))$\\
The units made of China-clay, glass, ceramic, mud, etc., are
kept in showrooms as heaped stocks. Due to this, units at
the bottom are under stress and as a result these units, break or get cracked. Thus breakability
rate of such an item depends on the accumulated stress of
stock level. So here  breakablility rate is dependent on stock variable x(t). So the breakability rate $B(x)$ is of the form
\begin{eqnarray}
B(x)=b_1 x^{\gamma} \,\,\, \mbox{and} \,\,\,\,\, b_1 \mbox{=constant  with} \,\,\,\,\, \gamma>0; 
\end{eqnarray}
 \noindent The unit
production cost is considered as a function of produced-quantity,
raw material cost, wear-tear and development costs(cf. Khouja [\ref{Khouja}]). So the total
production cost  is
\begin{eqnarray}
 c_{u} (t)u(t)=\bigg(c_{10}u(t) +{c_{d} (l
)} +\beta _{10} u^2 (t)\bigg)\end{eqnarray} \noindent
\noindent  The total holding cost over the finite time interval
$[0,\,T]$
 for the stock $x(t)$ is
\begin{eqnarray}
 \int _{0}^{T} h(t) x (t)\,dt \nonumber\ 
 \end{eqnarray}
 where $ h(t)=a+b t^n $ and $ a,b $ are being constants. \\
For a production system set-up cost is taken as normally constant. But, if dynamic
production rate is considered, some machineries, etc., are to be
set-up and maintained in such a way that the production
system can stand with the pressure of increasing demand.
Thus  a part of set up cost per unit time is linearly proportional to production rate and
hence the form is of the form   $\frac{s_1 +s_2 u(t)}{T}$. Here $s_1$ and $ s_2$  are constants.
Therefore the set up cost $s(t)$ for the model is taken as
\begin{eqnarray}
s(t)=\frac{s_1 +s_2 u(t)}{T}
\end{eqnarray}
 Let $p$  be the selling price of an unit.  Then the revenue from market demand is 
$p d(t)$. Thus the problem reduces to maximization of the profit function $J$
subject to the  constraint satisfying the dynamic
production-demand relation.
\begin{eqnarray}
&&\mbox{Max} J=\int _{0}^{T}\bigg(p d(t) - h(t) x(t)-c_{u }
 (t) u (t)-s(t)\bigg) dt\nonumber\\
&&\mbox{sub\, to}\,\,
 \frac{d}{dt}\bigg(x(t)\bigg)= u (t)-d(t)-B(x)\label{5}\\
  &&\mbox{ where \,} u(t)\ge 0\,\,  ,x(t)\ge0\label{6}\\
\nonumber\end{eqnarray}
\noindent Then the expression (\ref{5}) is written as
\begin{eqnarray}
 J_{max}\bigg({x}(t),\dot{x}(t),t\bigg)
 &=&\int _{0}^{T}\bigg(pd(t)-h(t)\, x(t)-\bigg(c_{10}+\frac{ s_2}{T}\bigg)
 \bigg({\dot{x}(t)+d(t)+B(x)}\bigg) \nonumber\\&-&{ \beta_{10}}\bigg({\dot{x}(t)+d(t)+B(x)}\bigg)^2 -\bigg({c_{d}(l)}+\frac{ s_1}{T}\bigg)\bigg) dt   \,\, \mbox{ where \,} \,  x(t)\ge0\label{7}
 \end{eqnarray}
 \noindent The above problem (\ref{7}) is defined as an optimal
control problem with state variable $x(t)$. Here, (\ref{7}) contains $u(t)$ implicitly.\\ 
\noindent Using Euler-Lagrange's equation, Fixed-Final Time and Fixed State System (i.e, here final
time $T$
is specified and $x(t)\equiv x, \dot x(t)\equiv \dot x $, $x(0)=0$ , $ x(T)=0$); we
have
\begin{eqnarray}\frac{\partial
J}{\partial {x}}&-&\frac{d\ }{dt}\bigg(\frac{\partial
J}{\partial \dot{x}}\bigg)=0\label{8}\\
\nonumber\end{eqnarray}  
\noindent Using equation (\ref{8}) in (\ref{7}), we get \\
\begin{eqnarray}
-h(t)-\bigg(c_{10}+\frac{s_2}{T}\bigg)b_1 \gamma x^{\gamma-1}-2 \beta_{10}\bigg(\dot x(t)+d(t)+b_1 x^{\gamma}\bigg)b_1 \gamma x^{\gamma-1}+2 \beta_{10}\frac{d}{dt}\bigg(\dot x(t)+d(t)+b_1 x^{\gamma}\bigg)=0\end{eqnarray}

Representing $x(t)\equiv x$,
 \begin{eqnarray}
\frac{d^2 }{dt^2}\bigg(x\bigg)-b^2_1 \gamma x^{2 \gamma-1}-\bigg(c_{10}+\frac{s_2}{T}+2 \beta_{10} d(t)\bigg)b_1 \gamma x^{\gamma-1} =\frac {1}{2\beta_{10}}\bigg(h(t)-{2\beta_{10}}(d_2+2 d_3 t)\bigg)
\end{eqnarray}
This is a second order first degree differential equation of the independent variable t.\\

\noindent Therefore, solving the above we can get values of 
\begin{eqnarray}
x(t)=\mbox{Complementary Function + Particular Integral }
\end{eqnarray}
 Here the initial condition  $x(0)=0$  and final condition
  $x(T)=0$. Using these conditions, we can get the constants of integration. From(1), we have the production rate 
 \begin{eqnarray}
 u(t)=\frac{d}{dt}\bigg(x(t)\bigg) +d_1 +d_2 t+d_3 t^2 +b_1 x(t)  \end{eqnarray}

  \noindent Then the  profit functional $J$ takes the form $J_{max}$
 \begin{eqnarray}
  &=&\int _{0}^{T}\bigg(p(d_1 +d_2 t+d_3 t^2)- h(t) \, x(t)-c_{10} u(t)-c_d(l)-\beta_{10} u^2(t)-\frac{s_1}{T}-\frac{s_2}{T} u(t)\bigg) dt\nonumber\\
 &=&\int _{0}^{T}\bigg(p(d_1 +d_2 t+d_3 t^2)- h(t) \, x(t)-\bigg(c_{10}+\frac{s_2}{T}\bigg) u(t)-\beta_{10} u^2(t)-\bigg(c_d(l)+\frac{s_1}{T}\bigg)\bigg) dt
\end{eqnarray}

\subsection{Models with particular values of $n\,\,\, \mbox{and}\,\,\, \gamma $:} 
\subsubsection{Model-1a: Formulation of optimal control  models for dynamic demand with holding cost $h(t)=a+bt$ for $n=1$ and linear breakability $B(x,t)=b_1 x(t) $ for $\gamma=1$ }
\indent In this model,  the differential equation for stock level
{$x(t)$} regarding above system
during a fixed time-horizon, T is
\begin{eqnarray} \frac{d}{dt}\bigg(x(t)\bigg)= u (t)-d(t)-b_1 x(t) \,\,\,\,\,\,\mbox{where} \,\,\,\,\, d(t)=d_1 + d_2 t+d_3 t^2 \label{14}
\end{eqnarray}
Proceeding as before, 
 we have 
\begin{eqnarray}
 J_{max}\bigg({x}(t),\dot{x}(t),t\bigg)
 &=&\int _{0}^{T}\bigg(pd(t)-h(t)\, x(t)-\bigg(c_{10}+\frac{ s_2}{T}\bigg)\bigg({\dot{x}(t)+d(t)+ b_1 x(t)}\bigg) \nonumber\\
 &-&{ \beta_{10}}\bigg({\dot{x}(t)+d(t)+ b_1 x(t)}\bigg)^2 -\bigg({c_{d}(l)}+\frac{ s_1}{T}\bigg)\bigg) dt   \,\, \mbox{ where \,} \,  x(t)\ge0 \label{15}
 \end{eqnarray}
 \noindent The above problem (\ref{15}) is defined as an optimal
control problem with state variable $x(t)$. Here, (16) contains $u(t)$ implicitly.\\ 
\noindent Using Euler-Lagrange's equation, Fixed-Final Time and Fixed State System (i.e, here final
time $T$
is specified and $x(t)\equiv x, \dot x(t)\equiv \dot x $, $x(0)=0$ , $ x(T)=0$); we
have
\begin{eqnarray}\frac{\partial
J}{\partial {x}}&-&\frac{d\ }{dt}\bigg(\frac{\partial
J}{\partial \dot{x}}\bigg)=0\label{16}\\
\nonumber\end{eqnarray}  
\noindent Using equation (\ref{16}) in (\ref{15}), we get \\
\begin{eqnarray}
-h(t)-\bigg(c_{10}+\frac{s_2}{T}\bigg)b_1-2 \beta_{10}\bigg(\dot x(t)+d(t)+b_1 x(t)\bigg)b_1
+2 \beta_{10}\frac{d}{dt}\bigg(\dot x(t)+d(t)+b_1 x(t)\bigg)=0\end{eqnarray}
It can be written as
\begin{eqnarray}
\frac{d^2 }{dt^2}\bigg(x(t)\bigg)-b^2_1  x(t)= \frac{1}{2 \beta_{10}}\bigg(a_{11} +a_{22} t + a_{33} t^2\bigg)
\end{eqnarray}
 where $a_{11}=(a+b_1 c_{10}+ b_1 \frac{s_2}{T} +2 \beta_{10} b_1 d_1-2\beta_{10} d_2); \,\,\,
a_{22}=(b-4 \beta_{10} d_3 +2 \beta_{10} b_1 d_2);   \,\,\,
a_{33}=2 \beta_{10} b_1 d_3 $\\
Then we have
\begin{eqnarray}
\frac{d^2 }{dt^2}\bigg(x(t)\bigg)-b^2_1  x(t)=f(t)  \,\,\,\,\mbox{where}\,\,\,f(t)=\frac{1}{2 \beta_{10}}\bigg(a_{11} +a_{22} t + a_{33} t^2\bigg)
\end{eqnarray}
This is a second order first degree differential equation of the independent variable t.\\
Therefore,
\begin{eqnarray}
x(t)&=& A_1 e^{b_1 t} +  B_1 e^{-b_1 t}+ \frac{f(t)}{D^2-b^2_1 };  D\equiv \frac{d}{dt}\nonumber\\
&=& A_1 e^{b_1 t} +  B_1 e^{-b_1 t}-\frac{f(t)}{b^2_1}-\frac{a_{33}}{b^4_1 \beta_{10}}\label{21}
 \end{eqnarray}
 where $A_1$ and $B_1$ are constants. 
 Therefore  $A_1$ and $B_1$  is given by the following relation. \\
  $A_1 +  B_1 = \frac{f(0)}{b^2_1}+\frac{a_{33}}{b^4_1 \beta_{10}}
$ and 
$ A_1 e^{b_1 T} +  B_1 e^{-b_1 T}=\frac{f(T)}{b^2_1}+\frac{a_{33}}{b^4_1 \beta_{10}}
$.\\
 From(\ref{14}), we have the production rate 
 \begin{eqnarray}
 u(t)= A_1 b_1 e^{b_1 t}-B_1  b_1 e^{-b_1 t} -\frac{a_{22}+2 a_{33}t }{2  b^2_1  \beta_{10}} +d_1 +d_2 t+d_3 t^2 +b_1 x(t) \end{eqnarray}
 It can be written as
  \begin{eqnarray}
 u(t)&=& 2 A_1 b_1 e^{b_1 t} -\bigg(\frac{a_{22} }{2  b^2_1  \beta_{10}}+\frac{a_{11}}{2 b_1 \beta_{10}}+\frac{a_{33}}{b^3_1 \beta_{10}}-d_1\bigg) -\bigg(\frac{ a_{33}}{  b^2_1  \beta_{10}} +\frac{a_{22}}{2 b_1 \beta_{10}}-d_2\bigg) t \label{22}
 \end{eqnarray}

  Then the  profit functional $J$ takes the form $J_{max}$
 \begin{eqnarray}
  &=&\int _{0}^{T}\bigg(p(d_1 +d_2 t+d_3 t^2)- h(t) \, x(t)-c_{10} u(t)-c_d(l)-\beta_{10} u^2(t)-\frac{s_1}{T}-\frac{s_2}{T} u(t)\bigg) dt\nonumber\\
 &=&\int _{0}^{T}\bigg(p(d_1 +d_2 t+d_3 t^2)- h(t) \, x(t)-\bigg(c_{10}+\frac{s_2}{T}\bigg) u(t)-\beta_{10} u^2(t)-\bigg(c_d(l)+\frac{s_1}{T}\bigg)\bigg) dt\nonumber\\
  &=&p(d_1 T +d_2 \frac{T^2}{2}+d_3 \frac{T^3}{3})-\bigg(c_d(l)+\frac{s_1}{T}\bigg)T\nonumber\\
  &-&a\bigg(\frac{A_1}{b_1} (e^{b_1 T}-1)+\frac{B_1}{b_1} (1-e^{-b_1 T})-\frac{1}{2 {b^2_1} \beta_{10}}\bigg(a_{11}T +a_{22} \frac{T^2}{2} + a_{33} \frac{T^3}{3}\bigg)-\frac{a_{33}}{b^4_1 \beta_{10}} T\bigg)+\frac{a_{33} b T^2}{b^4_1 2 \beta_{10}}\nonumber\\
  &-&b\bigg( A_1(\frac{T}{b_1} e^{b_1 T}+\frac{1}{b^2_1}(1-e^{b_1 T})) +  B_1 (-\frac{T}{b_1} e^{-b_1 T}+\frac{1}{b^2_1}(1-e^{-b_1 T}))-\frac{1}{2 {b^2_1} \beta_{10}}\bigg(a_{11}\frac{T^2}{2} +a_{22} \frac{T^3}{3} + a_{33} \frac{T^4}{4}\bigg)\bigg)\nonumber\\
  &-&C_{11}\bigg(2 A_1  (e^{b_1 T}-1) -\bigg(\frac{a_{22} }{2  b^2_1  \beta_{10}}+\frac{a_{11}}{2 b_1 \beta_{10}}+\frac{a_{33}}{b^3_1 \beta_{10}}-d_1\bigg)T -\bigg(\frac{ a_{33}}{  b^2_1  \beta_{10}} +\frac{a_{22}}{2 b_1 \beta_{10}}-d_2\bigg) \frac{T^2}{2}\bigg)\nonumber\\
 &-& \beta_{10} \bigg(M^2_1 T +M^2_2 \frac{T^2}{2}+M^2_3  \frac{(e^{2 b_1 T}-1)}{2b_1}+ M_1 M_2 T^2-2M_2 M_3 (\frac{T}{b_1} e^{b_1 T}+\frac{1}{b^2_1}(1-e^{b_1 T}))-2M_1 M_3 \frac{(e^{ b_1 T}-1)}{b_1} \bigg)\nonumber\\ \label{23}
 \end{eqnarray}
 \begin{eqnarray}
  &&\mbox{Where}\nonumber\\
  && M_3=-2 A_1 b_1 e^{b_1 T} ; M_1=\frac{a_{22} }{2  b^2_1  \beta_{10}}+\frac{a_{11}}{2 b_1 \beta_{10}}+\frac{a_{33}}{b^3_1 \beta_{10}}-d_1;
  M_2=\frac{ a_{33}}{  b^2_1  \beta_{10}} +\frac{a_{22}}{2 b_1 \beta_{10}}-d_2 \nonumber\\ 
  &&C_{11}=c_{10}+\frac{s_2}{T}; 
 \nonumber\\ \label{24}
\end{eqnarray}
\subsubsection{Model-1b: Formulation of optimal control  models for dynamic demand with holding cost $h(t)=a+bt$ with linear breakability $B(x,t)=b_1 x(t) $, $ b_1=0$}
If we take $b_1= 0$ in the
above model (Model-1a), then the present model represents an
inventory model without deterioration. As $b_1$ appears in the
denominator of the expression (\ref{21}), (\ref{22}), (\ref{23}) and (\ref{24}), it is not possible to put
$b_1$ directly  to obtain that expressions and the profit expression also. Thus, in
this case, the corresponding profit function is calculated omitting
$B(x,t)=b_1 x(t) $  and then proceeding in the same way as in Model-1a.
\begin{eqnarray}
a+bt- 2 \beta_{10} \frac{d\ }{dt} \bigg( \dot{x}(t)+d(t)\bigg)=0
\end{eqnarray} 
Representing $x(t)\equiv x$,
 \begin{eqnarray}
\frac{d^2 x}{dt^2}  =f(t)  \,\,\,\,\mbox{where}\,\,\, f(t)=-2 d_3 t+\frac{ b t}{2\beta_{10}}  +\frac{ a}{2\beta_{10}}-d_2
\end{eqnarray}
This is the second order first degree differential equation of the independent variable t.\\
\noindent Therefore, solving the above we can get  
\begin{eqnarray}
x(t)&=&A+Bt -(2d_3 -\frac{ b }{2\beta_{10}})\frac{t^3}{6}+(\frac{ a}{2\beta_{10}}-d_2) \frac{t^2}{2}
\end{eqnarray}
 Where $A$ and $B$ are constants. Here the initial condition  $x(0)=0$  and final condition
  $x(T)=0$. Using these conditions, we can get the values of  $A=0$  and $B=(2d_3 -\frac{ b }{2\beta_{10}})\frac{(T)^2}{6}+(\frac{ a}{2\beta_{10}}-d_2) \frac{T}{2}
$. \\
Then, we have the production rate 
 \begin{eqnarray}
 u(t)=B+ d_1 +\frac{ a}{2\beta_{10}} t+\frac{ b }{2\beta_{10}}\frac{t^2}{2}  
   \end{eqnarray}
\noindent Then the  profit functional $J$ takes the form $J_{max}$
 \begin{eqnarray}
 &=&\int _{0}^{T}\bigg(p(d_1 +d_2 t+d_3 t^2)- h(t) \, x(t)-\bigg(c_{10}+\frac{s_2}{T}\bigg) u(t)-\beta_{10} u^2(t)-\bigg(c_d(l)+\frac{s_1}{T}\bigg)\bigg) dt\nonumber\\
 &=&p(d_1 T +d_2 \frac{T^2}{2}+d_3 \frac{T^3}{3})-\bigg(c_d(l)+\frac{s_1}{T}\bigg)T-a\bigg( A T+B \frac{T^2}{2} -(2d_3 -\frac{ b }{2\beta_{10}})\frac{T^4}{24}+(\frac{ a}{2\beta_{10}}-d_2) \frac{T^3}{6}
\bigg)\nonumber\\
 &-&b\bigg( A \frac{T}{2}+B \frac{T^3}{3} -(2d_3 -\frac{ b }{2\beta_{10}})\frac{T^5}{30}+(\frac{ a}{2\beta_{10}}-d_2) \frac{T^4}{8}\bigg)-\bigg(c_{10}+\frac{s_2}{T}\bigg)\bigg( B T+ d_1 T +\frac{ b }{2\beta_{10}}\frac{T^3}{6}+\frac{ a}{2\beta_{10}} \frac{T^2}{2}\bigg)\nonumber\\ 
 &-&\beta_{10}\bigg( (B+ d_1)^2 +(\frac{ b }{2\beta_{10}})^2\frac{T^5}{20}+(\frac{ a}{2\beta_{10}})^2 \frac{T^3}{3}+\frac{ b(B+d_1) }{\beta_{10}}\frac{T^3}{6}+\frac{ a(B+d_1)}{2 \beta_{10}} T^2 +\frac{ a b}{4\beta^2_{10}} \frac{T^4}{8}\bigg)  
\end{eqnarray}
\newpage
\section{Numerical Experiments}
\noindent To illustrate the models, we assume the following input data for Models -1a and 1b. For Model-1a, $b_1$ has taken as $0.02$. and the other inputs are given in the following Table-\ref{t1} 
\subsection{Input:}
 Inputs are given in the following Table-\ref{t1} 
\begin{table}[h]
\begin{center} \caption{Input Values}
\scriptsize
\begin{tabular}{|cc c  c ccc c |} \hline
$L$   &$ N$&$  cd(l)$ & $ c_{10}$& $ \beta_{10} $&$p \mbox{(in \$)}$& $ s_1$ &
\\\hline\label{t1}
$ 40$ & $60$    & $100$     
&$0.7$ & $0.5$& $200$&  $ 10$ &\\\hline
$a$   &$ b $&$ d_1$ & $ d_2$& $ d_3 $&$T$& $ s_2$ &
\\\hline
$ 3$ & $0.2$        & $7$ 
&$4$ & $2$& $12$&  $ 3$ &\\\hline
\end{tabular}
\end{center}
\end{table}
\subsection{Output:}
\begin{figure}
\centering
\includegraphics[width=7.5cm, height=6.5cm]{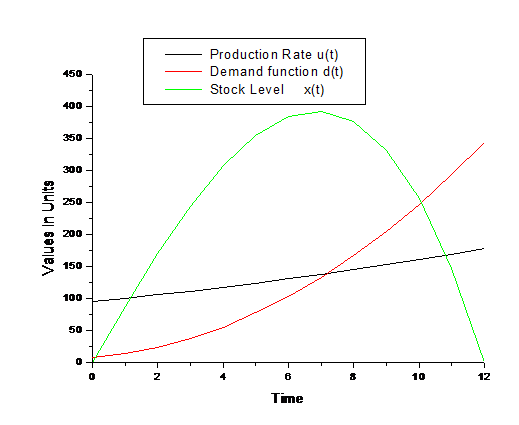}
\caption{Time Vs Production, Stock Level and Demand for Model-1a.}\label{f1}
\end{figure}

\begin{figure}
\centering
\includegraphics[width=7.5cm, height=6.5cm]{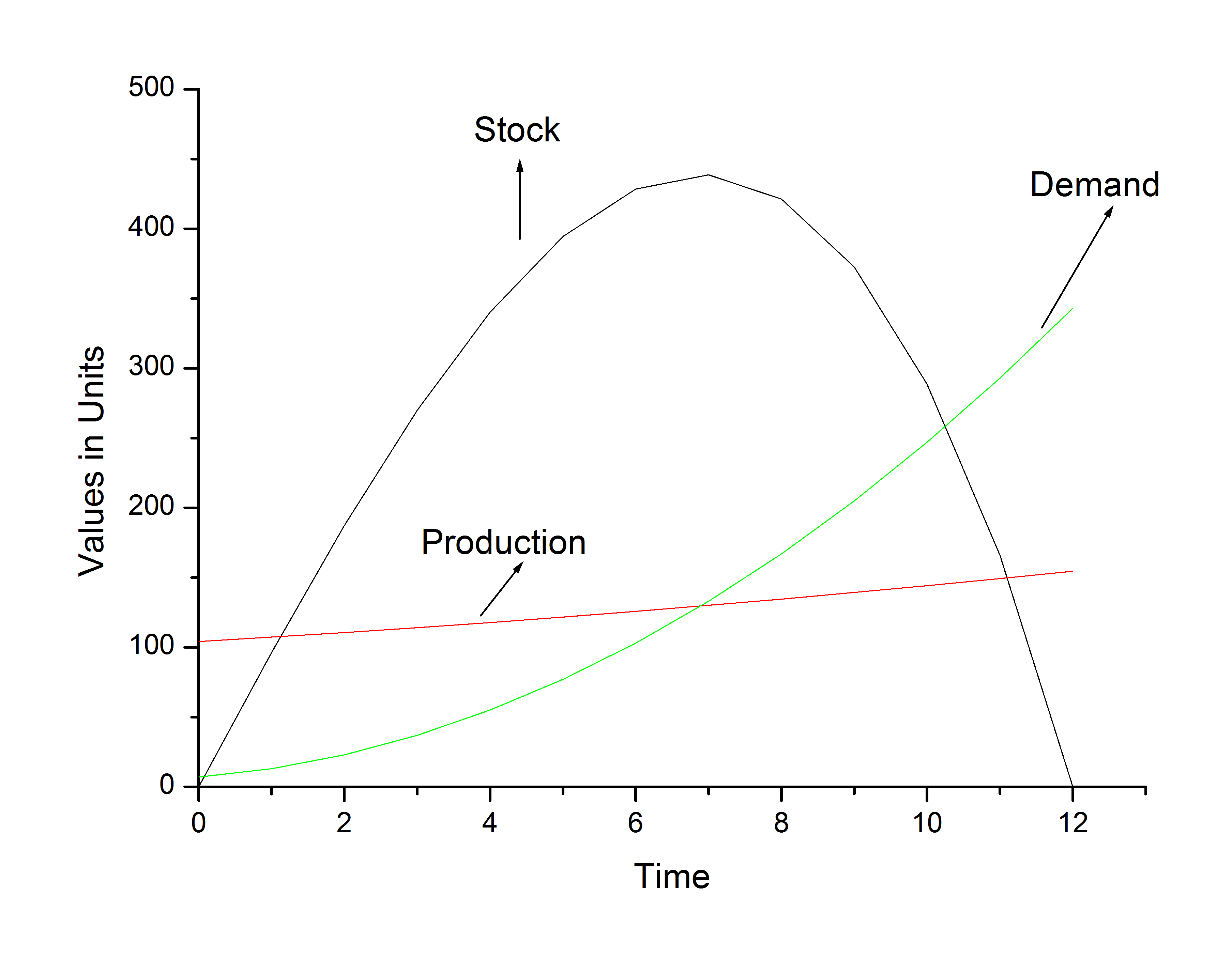}
\caption{Time Vs Production, Stock Level and Demand for Model-1b.}\label{f2}
\end{figure}

\noindent Using the above input data, the profit functions of Model$-1a$ and $1b$ are maximized and
with the help of GRG (LINGO-11.0), the unknowns, $x(t),u(t)$ and maximum profits etc. are
evaluated. Here  profit is $\$180913.30$ for Model-1a when breakability $b_1=0.02$ and profit is          $ \$ 247007.30$ for Model-1b when  breakability $b_1=0.00$(i.e the word breakability is not applicable). The numerical values of $ u(t), d(t), x(t)$
for the Model-1a (where  breakability is only dependent on stock and linear holding cost) and Model-$1b$ are given in Tables-\ref{t2} and -\ref{t3}. These are also graphically depicted in Figure-\ref{f1} and \ref{f2}. For the value of $b_1=0.11$, the unknowns $x(t),u(t)$ and $d(t)$ are evaluated which are given in Table-\ref{t4} and these are depicted in Figure-\ref{f3}. For the different values of breakability co-efficient $b_1$, profits are given in Table-\ref{t5}.
\begin{table}[h]
\begin{center}
\caption{Values of $u(t), d(t)$ and $x(t)$;  for Model-1a}
\resizebox{\linewidth}{!}{{\scriptsize 
\begin{tabular}{|c|c c c c c c c c c c c c c|} \hline
$t$   & $0$ & $1$ & $2$ & $3$ & $4$ &$5$ & $6$ & $7$& $8$& $9$ & $10$ & $11$ & $12$
\\\hline\label{t2}
$u(t)$ & $94.98$ & $100.05$ & $105.43$ &$111.16$ & $117.17$& $123.44$ &$130.09$ & $137.08$& $144.41$  & $152.10$& $160.14$ &$168.54$  & $177.32$ 
\\\hline
 $d(t)$ & $7.00$ & $13.00$ & $23.00$ & $37.00$ & $55.00$ & $77.00$  & $103.00$  & $133.00$  & $167.00$    & $205.00$  &$247.00$& $293.00$ & $343.00$ 
\\\hline
 $x(t)$ & $0$ & $86.95$&  $169.44$ &$243.86$ & $306.76$ & $354.70$ & $384.35$ & $392.44$ & $375.77$  & $331.22$& $255.72$ &$146.28$  & $0 $
\\\hline
\end{tabular}
}}
\end{center}
\end{table}
\begin{table}
\begin{center}\scriptsize 
\caption{Values of $u(t), d(t)$ and $x(t)$;  for Model-1b} 
\begin{tabular}{|c|c c c c c c c c c c c c c|} \hline
$t$   & $0$ & $1$ & $2$ & $3$ & $4$ &$5$ & $6$ & $7$& $8$& $9$ & $10$ & $11$ & $12$
\\\hline\label{t3}
$u(t)$ & $104.20$ & $107.30$ & $110.60$ &$114.10$ & $117.17$& $121.70$ &$125.80$ & $130.10$& $134.60$  & $139.30$& $144.20$ &$149.30$  & $154.60$ 
\\\hline
 $d(t)$ & $7.00$ & $13.00$ & $23.00$ & $37.00$ & $55.00$ & $77.00$  & $103.00$  & $133.00$  & $167.00$    & $205.00$  &$247.00$& $293.00$ & $343.00$ 
\\\hline
 $x(t)$ & $0$ & $96.06$&  $187.33$ &$270.00$ & $340.26$ & $394.33$ & $428.66$ &$438.66$& $421.33$ & $372.60$  & $288.66$& $165.732$&   $0 $
\\\hline
\end{tabular}
\end{center}
\end{table}
\section{Sensitivity Analysis}
\begin{figure}
\centering
\includegraphics[width=7.5cm, height=6.5cm]{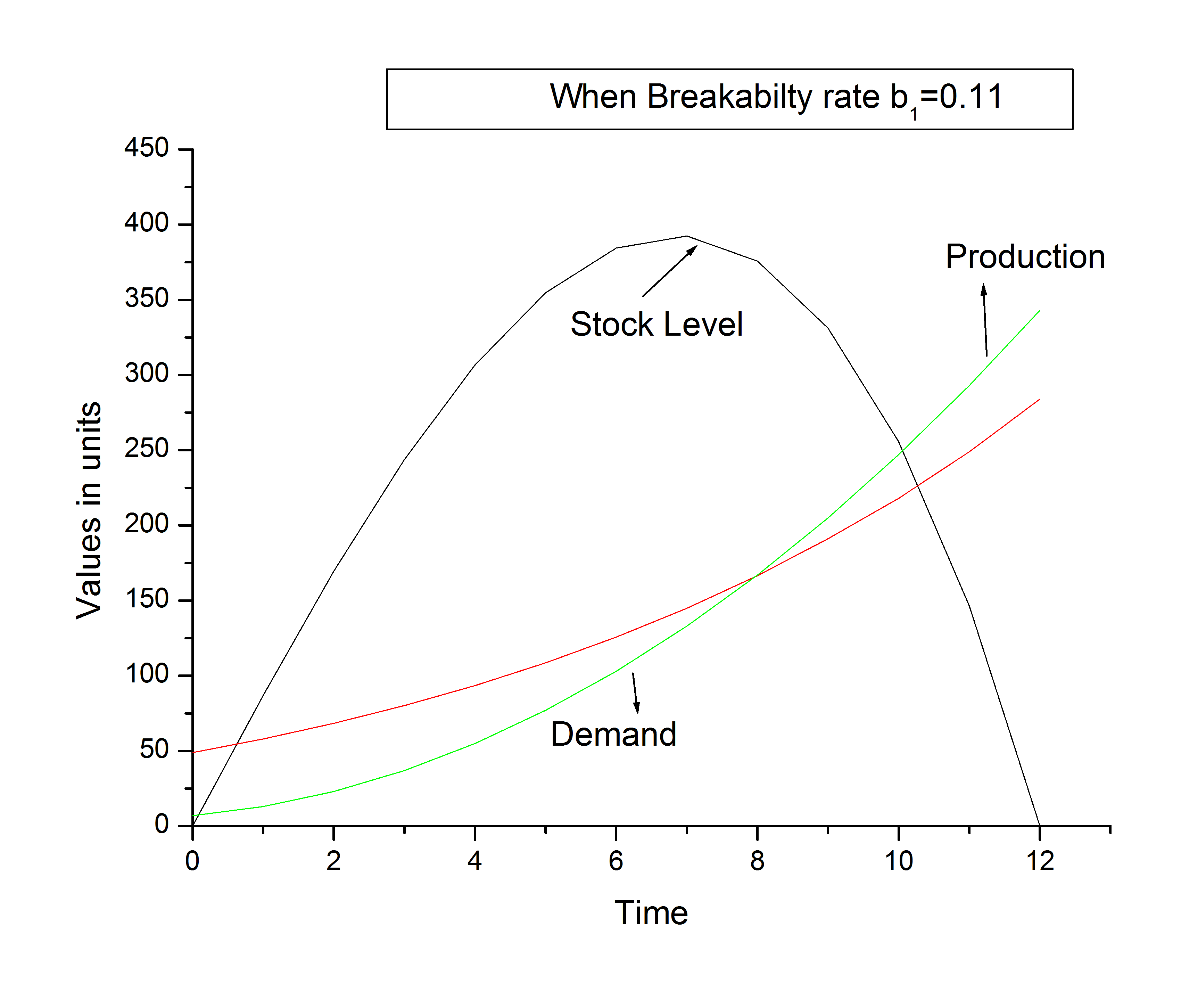}
\caption{Time Vs Production, Stock Level and Demand for Model-1a when breakability $b_1=0.11$ }\label{f3}
\end{figure}

\begin{figure}
\centering
\includegraphics[width=7 cm, height=6.5 cm]{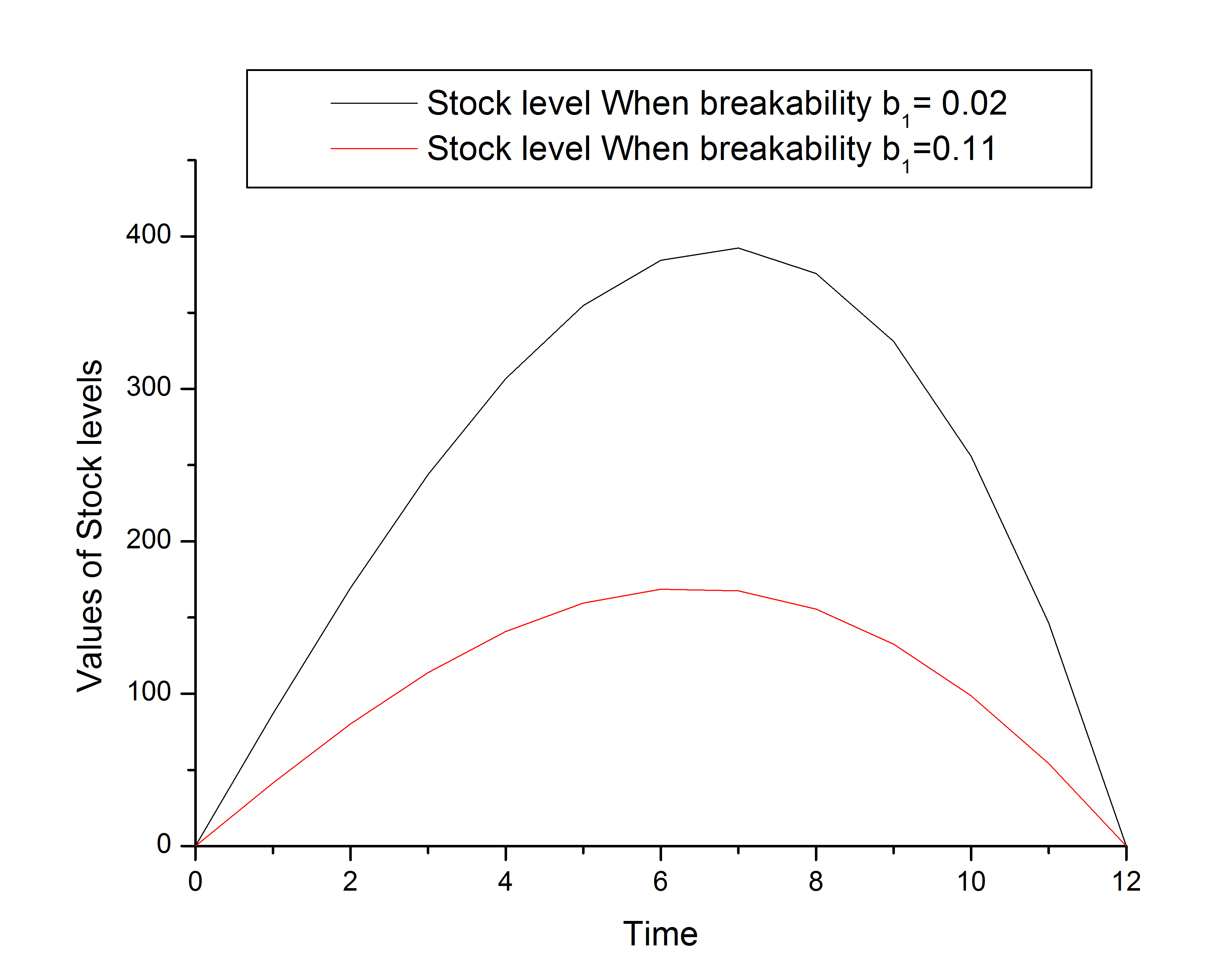}
\caption{Comparison of Time Vs Stock Level for Model-1a for the values of $b_1=0.11$  and $b_1=0.02$.}\label{f4}
\end{figure}

\begin{figure}
\centering
\includegraphics[width=7 cm, height=6.5 cm]{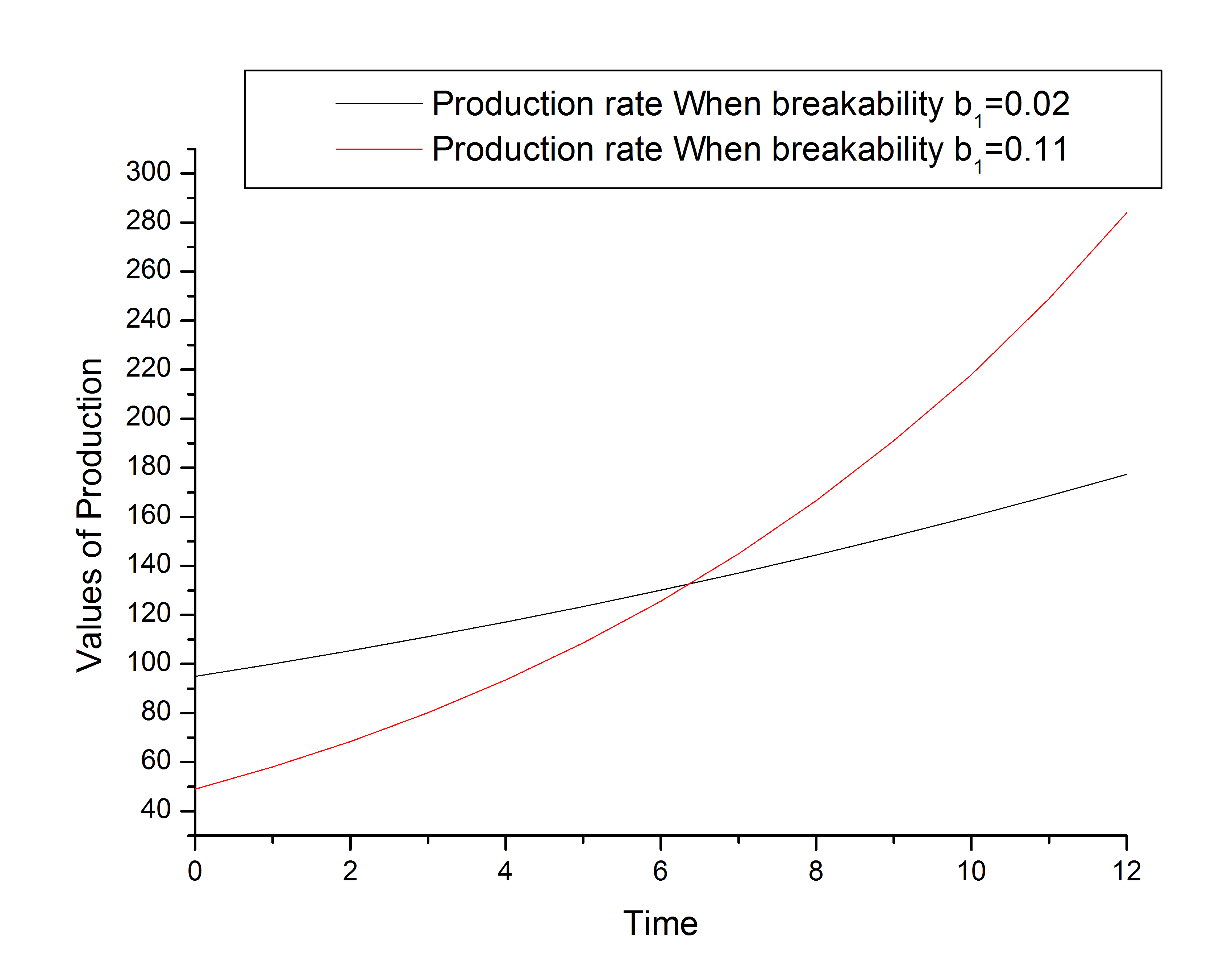}
\caption{Comparison of Time Vs Production rate for Model-1a for the values of $b_1=0.11$  and $b_1=0.02$.}\label{f5}
\end{figure}

For production system, a general belief that the system is more and more perfect as breakability of the system decreases.  From Table-\ref{t4}, it is seen that profit decreases as breakability rate increases. This is as per our expectation. The breakability constant $b_1$ has taken in the input as $0.02$. For this value of $b_1$, the profit for Model-1a is 180913.30 $\$$. With $b_1=0.11$ , the value of profit for Model-1a is 153447.70 $\$$. For that value of $b_1$, the unknowns, $x(t),u(t)$ and $d(t)$ are evaluated which are given in Table-\ref{t4} and these are depicted in figure-\ref{f3}. 
Taking the value of $b_1$ as $b_1=0.11$  and$b_1=0.02$ for Model-1a, the comparison  diagram for production and stock  has been drawn in figures-\ref{f4} and \ref{f5}.
\begin{table}
\begin{center}
\caption{Values of $u(t), d(t)$ and $x(t)$;  for Model-1a, [The value of $b_1$ as $b_1=0.11$ , the value of profit for Model-1a is 153447.7 $\$$]}
\resizebox{\linewidth}{!}{{\scriptsize
\begin{tabular}{|c|c c c c c c c c c c c c c|} \hline
$t$   & $0$ & $1$ & $2$ & $3$ & $4$ &$5$ & $6$ & $7$& $8$& $9$ & $10$ & $11$ & $12$
\\\hline\label{t4}
$u(t)$ & $48.96$ & $58.04$ & $68.38$ &$80.15$ & $93.48$& $108.58$ &$125.65$ & $144.92$& $166.63$  & $191.09$& $218.56$ &$249.45$  & $284.18$ 
\\\hline
 $d(t)$ & $7.00$ & $13.00$ & $23.00$ & $37.00$ & $55.00$ & $77.00$  & $103.00$  & $133.00$  & $167.00$    & $205.00$  &$247.00$& $293.00$ & $343.00$ 
\\\hline
 $x(t)$ & $0$ & $41.44$&  $80.15$ &$113.90$ & $140.83$ & $159.44$ & $168.59$ & $167.44$ & $155.48$  & $132.49$& $98.59$ &$54.18$  & $0 $
\\\hline
\end{tabular}
}}
\end{center}
\end{table}
\begin{table}
\begin{center}
 \caption{The values of $b_1$ Vs Profit of Model-1a }
 \resizebox{\linewidth}{!}{{\scriptsize
\begin{tabular}{|c|cc c c c c cc|} \hline
$b_1$  &$0.01$ & $0.02$ & $0.03$ & $0.04$ & $0.05$ & $0.06$ &$0.07$&$0.08$ 
\\\hline\label{t5}
$\mbox{Profit} \,\,\, J\,\,\, in \,\,\, $\$$$ &$185131.50$& $180913.30$ & $176871.90$ & $173036.20$ &$169431.00$ & $166076.60$& $162988.70$&$160178$ 
\\\hline
 \end{tabular}
 }}
\end{center}
\end{table}
\newpage
\section{Discussion}
 For production system, a general belief is that cost of the produced item increases as breakability
units of the system increase. Consequently, Min cost of the production system increases. i.e,  total profit of the production system decreases.  These are as per our expection. In the Figures-\ref{f1}and \ref{f2}\,\,  for the models, initially production is more than demand and hence stock is built up after satisfying the demand. As demand is here dynamic which  rapidly increases with time, after some time,  surpasses the production and then excess demand is satisfied from the stock. So, the deficit in demand occurs, so the stock starts to decrease at this point and gradually decreases up to the end of the cycle
and finally becomes zero. It is observed that  profit is $ \$ 247007.30$ for Model-1b without  breakability which is more than the profit  $\$180913.30$ for Model-1a when breakability $b_1=0.02$. When the value of $ b_1$  increases, the value of profit decreases. This is clear in Table-4. For particular value of $b_1$ as $b_1=0.11$ , the value of profit for Model-1a is 153447.7 $\$$. This result is very natural and  is as per  our expectation.
\section{Conclusion}
Here a production-inventory model for a breakable item is developed.  Here time dependent production rate fetches more profit for  time-dependent demand. From the present model it can be concluded that optimal control of production rate reduces holding cost as well as damageability which in turn increases profit separately for breakable/deteriorating items.  The models are also solved taking some of the inventory costs. The present models can be extended to the rough, fuzzy-rough, random, fuzzy-random environment taking constant part of holding cost, set-up cost, etc. as uncertain in nature. As the formulation and solution are quiet general, the results can be obtained for other forms of dynamic demand such as linearly increasing demand, ramp demand etc. The model can be extended to include multi-item fuzzy inventory problem with fuzzy space and budget constraints.


\end{document}